%% file: main.tex
\documentclass[12pt]{article}
\usepackage{amsmath,amssymb,amsthm,amsfonts}
\usepackage{arxiv}
\usepackage{hyperref}
\usepackage{nicefrac}
\usepackage{bbm}
\usepackage{mathrsfs}
\usepackage{algorithm}
\usepackage{algorithmic}
\usepackage{times}
\usepackage{graphicx}
\usepackage{subcaption}
\usepackage{url}
\usepackage{booktabs}
\usepackage{secdot}
\usepackage{multirow}
\usepackage{multicol}

\input{math_symbols}

\title{Data-driven optimal prediction with control}
\usepackage{times}

\author{
Aleksandr Katrutsa\thanks{Corresponding author} \\
Skoltech, Moscow, Russia\\
AIRI, Moscow, Russia\\
\url{amkatrutsa@gmail.com} 
\And
Ivan Oseledets \\
AIRI, Moscow, Russia\\
Skoltech, Moscow, Russia\\
\url{oseledets@airi.net} 
\And
Sergey Utyuzhnikov\\
The University of Manchester \\
Manchester, M13 9PL, UK \\
\url{s.utyuzhnikov@manchester.ac.uk}
}

\begin{document}

\maketitle

\begin{abstract}
    This study presents the extension of the data-driven optimal prediction approach to the dynamical system with control.
The optimal prediction is used to analyze dynamical systems in which the states consist of resolved and unresolved variables.
The latter variables can not be measured explicitly. 
They may have smaller amplitudes and affect the resolved variables that can be measured.
The optimal prediction approach recovers the averaged trajectories of the resolved variables by computing conditional expectations, while the distribution of the unresolved variables is assumed to be known.
We consider such dynamical systems and introduce their additional control functions.
To predict the targeted trajectories numerically, we develop a data-driven method based on the dynamic mode decomposition.
The proposed approach takes the \emph{measured} trajectories of the resolved variables, constructs an approximate linear operator from the Mori-Zwanzig decomposition, and reconstructs the \emph{averaged} trajectories of the same variables.  
It is demonstrated that the method is much faster than the Monte Carlo simulations and it provides a reliable prediction.
We experimentally confirm the efficacy of the proposed method for two Hamiltonian dynamical systems.
\end{abstract}

\section{Introduction}
The modeling and analysis of many physical phenomena is based on the dynamical system framework~\cite{katok1995introduction}.
This study considers dynamical systems such that state vectors consist of resolved and unresolved variables.
Typically, the resolved variables correspond to the measured quantities, and their amplitudes dominate the amplitudes of the unresolved variables.
However, the unresolved variables implicitly affect the resolved variables and must be considered during the analysis of dynamical system trajectories.
The standard approach to revealing this implicit effect is called the \emph{optimal prediction}~\cite{chorin2002optimal}.
The optimal prediction framework produces the averaged trajectories of the resolved variables via the conditional expectation over the randomized unresolved variables with a known distribution.
The key tool for deriving the analytical solution in the optimal prediction is the so-called Mori-Zwanzig decomposition~\cite{lin2021data,mori1965transport,zwanzig1973nonlinear}.
This decomposition provides the reduced-order model for the resolved variables in the exact form, with the explicit effect of unresolved variables.
However, the analytical computing of the conditional expectation can be performed only for a limited number of dynamical systems~\cite{bernstein2007optimal,chorin2003averaging}.
Therefore, data-driven methods for revealing averaged dynamics of the resolved variables become more popular~\cite{garcia2023physics,dal2020data,kadeethum2021framework}.
In~\cite{lin2021data,curtis2021dynamic}, the authors combine the optimal prediction with the dynamic mode decomposition (DMD)~\cite{schmid2010dynamic,kutz2016dynamic}. This makes the optimal prediction approach potentially suitable for the analysis of complex physical systems when only measured data are available. 

In the present paper, we extend the technique proposed in \cite{curtis2021dynamic} to the dynamical systems with controls. 
Note that there is also a DMD approach for such systems  (DMDc)~\cite{proctor2016dynamic}, which adapts the DMD framework to the dynamical systems with additional control functions.
However, we show below that this approach can not accurately process dynamical systems with both \emph{control and unresolved variables}.
We aim to combine DMDc with the optimal prediction. 
Next, we call our approach the \emph{optimal prediction with control} (OPc) and demonstrate in experiments that it generates the averaged trajectories of the resolved variables faster than the straightforward Monte Carlo simulations.
Moreover, the OPc method provides the quasi-optimal linear operator, which reveals the stability properties of the considered dynamical system on average.
The numerical experiments show that OPc tracks the trend of the ground-truth averaged dynamics of the resolved variables.

\paragraph{Related works}
The data-driven methods for analyzing PDEs have recently attracted much interest~\cite{rudy2017data,berg2019data,raissi2019physics}.
The classical Dynamic Mode Decomposition method (DMD)~\cite{schmid2010dynamic} provides a simple yet efficient tool for extracting the underlying dynamic modes from the given snapshots.
DMD has a lot of applications in different domains like predictive modeling in agriculture~\cite{stasenko2023dynamic}, climate modeling~\cite{falkena2019derivation,ghil2020physics,palmer2019stochastic,franzke2020structure}, fluid dynamics~\cite{gouasmi2016towards,parish2016reduced}, molecular dynamics simulations~\cite{wang2019implicit,klippenstein2021introducing,volgin2018coarse}, etc.
The authors suggest multiple modifications of DMD to extract additional structures of solution~\cite{chen2012variants,le2017higher,azencot2019consistent,jovanovic2014sparsity,williams2015data} and improve the computational efficiency~\cite{erichson2019randomized,degennaro2019scalable,askham2018variable}.
In addition, the combination of the Gaussian process framework with DMD is presented in~\cite{tsolovikos2024dynamic}.
The theoretical foundation can be obtained from the theory of the Koopman operator~\cite{koopman1931hamiltonian,arbabi2017ergodic}. 
The analysis of the Koopman operator spectrum shows that the DMD method provides the linear finite-dimensional operator that approximates the discrete part of the Koopman operator spectrum~\cite{mauroy2020koopman,arbabi2017ergodic}. 
The Koopman operator approach is extended to the case of the dynamical systems with control, namely the Koopman theory with inputs and control (KIC), in~\cite{proctor2018generalizing}. 
However, the aforementioned methods do not consider dynamical systems with unresolved variables in the state vectors.

The optimal prediction framework provides the tools for analyzing dynamical systems with unresolved variables.
Moreover, in~\cite{curtis2021dynamic}, the authors establish the connection between the optimal prediction and specific form of the Dynamic Mode Decomposition (DMD).
This approach is modified in ~\cite{katrutsa2023extension}, where the t-model of the optimal prediction is implemented.
At the same time, incorporating a control function in the optimal prediction approach is still not well-studied.
In this work, we fill this gap and present both the theoretical derivation and practical tests of the optimal prediction with control (OPc) method.

The rest of the paper is organized as follows. 
In the next section the Mori-Zwanzig formalizm is delivered to identify the resolved variables from the Liouville equations. 
Then, in Section~\ref{sec:opc} the data-driven optimal prediction method is modified to adapt it to the problems with control. 
The efficacy of the proposed method is tested in numerical experiments presented in Section~\ref{sec:experiments} on dynamical systems with a constant and linear control. 

\section{Mori-Zwanzig projection}
\label{sec:main}

This section briefly presents the Mori-Zwanzig formalism that is used for estimating the averaged trajectories.
This formalism further leads to the optimization problem for computing the locally optimal operator which is able to provide the expectation for the trajectories of dynamical systems with uncertainties.

Let $\bu(\bz, t)$ be the state vector of a dynamical system with control.
This vector changes with time and composes a trajectory in the state space.
The evolution of the state vector is governed by the Liouville equation with the additional control item
\begin{equation}
    \begin{aligned}
    & \frac{\partial\bu}{\partial t} = \mathcal{L}\bu + \calB \bv,\\
    & \bu(\bz, 0) = \bu_0,
    \label{eq::liouv_general}
\end{aligned}
\end{equation}
where $\mathcal{L} = \sum_{j=1}^n f_j(\bz) \frac{\partial }{\partial z_j}$, $\bu_0 \in \bbR^n$ is the initial value, $\bv$ and $\calB$ are a control vector and control operator, respectively.


Assume that
\begin{equation}
    \calB \bv = \sum_{j=1}^n g_j \frac{\partial }{\partial z_j}\bu,
    \label{eq:Bv}
\end{equation}
where $\mathbf{g}$ is a control function. Then, the problem \eqref{eq::liouv_general} is equivalent to the following Cauchy problem for the Langevin equation:
 \begin{equation}
 \begin{split}
     &\frac{d\mathbf{\Phi}}{dt} = \mathbf{f}(\mathbf{\Phi})+ \mathbf{g}, \quad \mathbf{\Phi} \in \bbR^n,\\
     &\mathbf{\Phi}(0) = \bz.
 \end{split}
 \label{eq::caushy_char}
 \end{equation}
As can be seen, the equation~\eqref{eq::caushy_char} is the characteristic equation for the equation~\eqref{eq::liouv_general} with~\eqref{eq:Bv}. 

The equivalence of \eqref{eq::liouv_general}-\eqref{eq:Bv} and \eqref{eq::caushy_char} directly follows from the observation: $\bu(\bz, t) = \mathbf{\Phi}(\bz, t)$ if $\bu_0 = \bz$. 
In this case,  at $t=0$ $\calB \bv = \sum_{j=1}^n g_j \frac{\partial }{\partial z_j}\bu_0=\bg$. 
On the other hand, along each characteristic $\bu$ does not change. 
Hence, $\calB \bv$ does not change either. Thus, $\calB \bv$ is equal to $\bg$ everewhere.

In this work, we consider dynamical systems in which the state vector consists of the resolved and unresolved variables.
Without the loss of generality, let us split the state vector elements into \emph{resolved variables} $\bx(t)$ and \emph{unresolved variables} $\by(t)$ that can not be measured explicitly. Thus, $\bu(t) = [\bx(t), \by(t)]$.

Although the unresolved variables are supposed to be small, these variables affect the trajectories of the resolved variables and must be considered in modeling.
For this purpose, consider the projectors $\calP$ and $\calQ$ onto the spaces of resolved and unresolved variables, respectively.
In particular, the projector $\calP$ can be defined as $\bx(t) = \calP \bu(t) = \mathbb{E} [\bu(t) | \bx(t)]$, where the probability distribution of unresolved variables $\by(t)$ is assumed to be known in the computing of this conditional expectation.
In addition, note that $\calQ \bu(t) = \by(t)$ and $\calP + \calQ = \calI$, which is the identity operator.

Now, we can derive the dynamics equation for the resolved variable $\bx(t)$ only.
Using the semi-group notation for~\eqref{eq::liouv_general} solution $\bu(t) = e^{\calL t}\bu_0 + \bs(t)$ (where $\bs(t)$ is a particular solution of the non-homogeneous equation) and the projector~$\calP$, one can derive the following equation:
\begin{equation}
\begin{split}
\frac{\partial(\mathcal{P}e^{\mathcal{L}t}\bu_0 + \calP \bs)}{\partial t} &= \calP \calL e^{\calL t}\bu_0 + \calP \calL \bs + \calP\calB \bv \\
&= \calP e^{\calL t} \calP \calL \bu_0 + \calP e^{\calL t}\calQ\calL \bu_0 + \calP \calL \calP \bs + \calP \calL \calQ \bs + \calP \calB \bv ,
\end{split}
\end{equation}
Here the first item corresponds to the pure dynamics of the resolved variables.
At the same time, thanks to the Dyson operator identity
$e^{\calL t} = e^{\calQ \calL t} + \int_{0}^t e^{\calL (t - \tau)}\calP \calL e^{\calQ \calL \tau} d\tau$, the second item can be further re-written in the following form
\begin{equation*}
        \calP e^{\calL t}\calQ\calL \bu_0 = \calP \by(t) + \int_{0}^t \calP e^{\calL (t - \tau)} \calM (\tau, \bu_0) d\tau,
\end{equation*}
where $\calM (\tau, \bu_0) = \calP \calL e^{\calQ \calL \tau} \calQ \calL \bu_0 = \calP \calL \bN(\tau, \bu_0)$.
Note that $\bN(\tau, \bu_0)$ corresponds to the dynamics in the space of unresolved variables that can be interpreted as noise~\cite{curtis2021dynamic}.

The subspace of resolved variables is orthogonal to the subspace of unresolved variables that means $\calP \by(t) \equiv 0$.
The remaining term also $\calP\calL\calQ \bs \equiv 0$ since $\calL\calQ \bs$ belongs to the space of unresolved variables. Hence, its projection onto the space of resolved variables is zero. Thus, we arrive at the following equation for the resolved variables:
\begin{equation}
    \frac{\partial \bx}{\partial t} = \calP \calL \bx + \int_0^{t} \calP e^{\calL (t - \tau)} \calM(\tau, \bu_0) d\tau + \calP \calB \bv,
    \label{eq::mz_1}
\end{equation}
where the second term corresponds to the impact of unresolved variables on the resolved variables, and the last term is the projection of the control vector onto the space of resolved variables.

To derive the complete set of equations that operates only with the resolved variables, we have to exclude the action of projector $\calQ$ in the definition of the memory kernel $\calM(\tau, \bu_0)$.
Starting from the definition of $\bN(\tau, \bu_0)$ and Dyson operator identity, one can derive the following expression
\begin{equation*}
    \begin{split}
        \bN(\tau, \bu_0) &= e^{\calQ \calL \tau} \calQ \calL \bu_0 = e^{\calL \tau} \calQ \calL \bu_0 - \int_0^{\tau} e^{\calL (\tau - s)}\calP \calL e^{\calQ \calL s} \calQ \calL \bu_0 ds \\
        &= e^{\calL \tau} \bN(0, \bu_0) - \int_0^{\tau} e^{\calL (\tau - s)} \calM (s, \bu_0)ds.
    \end{split}
\end{equation*}
After multiplying both sides by $\calP\calL$, the latter equation is transformed to the following Volterra equation:
\begin{equation}
    \calM(\tau, \bu_0) + \int_0^{\tau} \calP e^{\calL (\tau - s)} \calL \calM (s, \bu_0)ds = \calP e^{\calL \tau} \calL \bN(0, \bu_0).
    \label{eq::mz_2}
\end{equation}
Now, we are ready to introduce the proposed data-driven method named optimal prediction with control (OPc), which is based on the discretization of~\eqref{eq::mz_1} and~\eqref{eq::mz_2} and estimates the finite-dimensional approximation $\bA_{opc}$ of the projected operator $\calP\calL$.

\section{Optimal prediction with control}
\label{sec:opc}

To derive the OPc method, we discretize~\eqref{eq::mz_1} and~\eqref{eq::mz_2} and briefly describe the auxiliary transformations.
A more detailed description of the following transformations is presented in~\cite{katrutsa2023extension}.
Let us start from equation~\eqref{eq::mz_1}, where we use the first-order discretization of the left-hand side and the trapezoidal rule for the integral item.
As a result, we obtain the following equation
\begin{equation}
\begin{split}
     \frac{\bx_{n+1} - \bx_n}{\Delta \tau} &= \calP\calL \bx_n + \Delta \tau \left( \sum_{k=1}^{n-1} \calP e^{(n - k)\Delta\tau\calL} \calM(k\Delta\tau, \bu_0) + \right.\\ &\left.\frac{1}{2}[\calP \calM(n\Delta\tau, \bu_0) + \calP e^{\calL n\tau}\calM(0, \bu_0)] \right) + \bw_n,
\end{split}
    \label{eq::mz_discr_1}
\end{equation}
where $\bx_n = \bx(n\Delta \tau)$ and $\bw_n = \calP \calB \bv(n \Delta \tau)$ are the snapshots of the resolved variables and control vectors, respectively.

Here, we assume that the operator $\calB$ is known, and we can operate with the transformed control vectors $\calP\calB\bv_n$ explicitly.
Assume that the target operator $\calP \calL$ can be approximated with the finite-dimensional operator $\bA_{opc} \approx \calP \calL $ that admits the eigendecomposition $\bA_{opc} = \bV \bLambda \bV^{-1}$.
Next, let us re-write~\eqref{eq::mz_discr_1} in the following form
\begin{equation}
    \hat{\bx}_{n+1} = (\bI + \Delta \tau\bLambda)\hat{\bx}_n + \frac{\Delta \tau^2}{2} \left(\widehat{\calM}_n + e^{n\Delta\tau\bLambda}\widehat{\calM}_0 + 2 \sum_{k=1}^{n-1}e^{(n-k)\Delta \tau \bLambda}\widehat{\calM}_k \right) + \hat{\bw}_n \Delta \tau,
    \label{eq::x_discr}
\end{equation}
where $\hat{\bx}_n = \bV^{-1}\bx_n$, $\widehat{\calM}_k = \bV^{-1}\calM_k$ and $\hat{\bw}_n = \bV^{-1}\bw_n$.

Now, we derive the discretized version of~\eqref{eq::mz_2} via the same trapezoidal rule and re-write it in the basis $\bV$:
\begin{equation}
    \widehat{\calM}_n = \left( \bI + \frac{\Delta\tau}{2}\bLambda \right)^{-1} \left( e^{\bLambda n\Delta\tau} \left(\bI - \frac{\Delta\tau}{2}\bLambda \right)\widehat{\calM}_0 - \Delta\tau \bLambda  \sum_{k=1}^{n-1} e^{\bLambda (n-k)\Delta\tau} \widehat{\calM}_k \right).
    \label{eq::mz_memory_discr}
\end{equation}
This equation admits the following recurrent form:
\[
\widehat{\calM}_n = e^{\Delta\tau\bLambda}\bM(\bLambda)\widehat{\calM}_{n-1},\; n \geq 1,
\]
where $\bM(\bLambda) = \bI - \Delta \tau \bLambda \left( \bI + \frac{\Delta\tau}{2}\bLambda\right)^{-1}$.
Thus, Eq.~\eqref{eq::x_discr} is transformed to the following form
\begin{equation}
        \hat{\bx}_{n+1} = (\bI + \Delta \tau \bLambda)\hat{\bx}_n + \frac{\Delta \tau^2}{2}e^{n\Delta\tau \bLambda}\left(\bM^n(\bLambda) + \bI + 2\sum_{k=1}^{n-1} \bM^k(\bLambda) \right)\widehat{\calM}_0 + \hat{\bw}_n \Delta \tau.
        \label{eq::x_discr_rec}
\end{equation}
Note that the expression with $\bM(\bLambda)$ in~\eqref{eq::x_discr_rec} can be simplified to $(\bM^n(\bLambda) - \bI)\left(-\frac{2}{\Delta\tau} \bLambda^{-1} \right)$.

Thus, after substitution of this expression in~\eqref{eq::x_discr_rec}, we arrive at the following equation:
\begin{equation}
    \hat{\bx}_{n+1} = (\bI + \Delta \tau \bLambda)\hat{\bx}_n -  \Delta \tau \bLambda^{-1} e^{n\Delta\tau \bLambda}(\bM^n(\bLambda) - \bI) \widehat{\calM}_0 + \hat{\bw}_n \Delta \tau.
    \label{eq::x_discr_lamb}
\end{equation}
Introduce the following notations: $\bar{\bLambda} = \bI + \Delta \tau \bLambda$ and $\bA = \bV\bar{\bLambda}\bV^{-1}$.
Then, we can re-write~\eqref{eq::x_discr_lamb} in the terms of $\bar{\bLambda}$ based on equality $\bLambda = \frac{\bar{\bLambda} - \bI}{\Delta\tau}$:
\begin{equation}
    \hat{\bx}_{n+1} = \bar{\bLambda}\hat{\bx}_n -  \Delta \tau^2 (\bar{\bLambda} - \bI)^{-1}e^{n(\bar{\bLambda} - \bI)}(\bM^n(\bar{\bLambda}) - \bI) \widehat{\calM}_0 + \hat{\bw}_n \Delta \tau,
    \label{eq::x_discr_lamb_hat}
\end{equation}
where $\bM(\bar{\bLambda}) = \bI - 2 (\bar{\bLambda} - \bI)(\bar{\bLambda} + \bI)^{-1}$ or in terms of the transition operator~$\bA$:
\begin{equation}
    \bx_{n+1} = \bA\bx_n -  \Delta \tau^2(\bA - \bI)^{-1} e^{n(\bA - \bI)}(\bM^n(\bA) - \bI) \calM_0 + \bw_n\Delta \tau,
    \label{eq::opc_generation}
\end{equation}
where $\bM(\bA) = \bI - 2 (\bA - \bI)(\bA + \bI)^{-1}$.

Now, we have the explicit relation between the two sequential snapshots $\bx_{n+1}$ and $\bx_n$, which leads to the following optimization problem to obtain the operator $\bA_{opc}$:
\begin{equation}
\bA_{opc} = \argmin_{\bA} \| \bY_+ - \bA \bX_- + \Delta \tau^2 \widetilde{\bM}(\bA, \bn) \|_F^2,
    \label{eq::mzdmd_final_opt}
\end{equation}
where $\bX_- = [\bx_{m-1}, \ldots, \bx_1]$ is the stacked snapshots and $\bY_+ = [\bx_m - \bw_{m-1}\Delta \tau, \ldots, \bx_2 - \bw_1\Delta \tau]$ is the stacked snapshots corrected by the control vectors. 
Both matrices $\bX_-$ and $\bY_+$ are of size $n \times (m-1)$.
In addition, denote the stacked vectors of corrections in~\eqref{eq::opc_generation} by $\widetilde{\bM}(\bA, \bn) = (\bA - \bI)^{-1}\bF(\bA, \bn)$ and 
\[
\bF(\bA, \bn) = 
\begin{bmatrix} \mathbf{0} & F_1(\bA, \bn) & \ldots & F_{m-2}(\bA, \bn) \end{bmatrix}
\]
is a matrix of size $n\times (m-1)$ such that $F_j(\bA, \bn) = e^{j(\bA - \bI)}(\bM^j(\bA) - \bI)\bn$. 
Here, we also introduce the vector~$\bn$ corresponding to a particular initialization of the memory term $\calM_0$.

Note that if one eliminates the correction term $\Delta \tau^2 \widetilde{\bM}(\bA, \bn)$ from~\eqref{eq::mzdmd_final_opt}, then the solution of the resulting optimization problem coincides with the solution of DMD with control (DMDc) method proposed in~\cite{proctor2016dynamic}.

\paragraph{If the operator $\calB$ is unknown.}
In this case, we can not include control vectors $\bw_n$ explicitly in the composing of the matrix $\bY_+$.
Instead, we can access the input vectors $\bv_n = \bv(n\Delta \tau)$ and should keep the unknown matrix $\bB_{opc} \approx \calP\calB$ in our derivation.
In the final step, this unknown matrix becomes one more variable in the resulting optimization problem in addition to the matrix $\bA$.
Thus, the adjustment of the aforementioned derivation procedure for unknown matrix $\bB$ leads to the following optimization problem
\begin{equation}
    \bA_{opc}, \bB_{opc} = \argmin_{\bA, \bB} \| \bX_+ - \bA \bX_- + \Delta \tau^2 \widetilde{\bM}(\bA, \bn) - \bB \bV_c \|_F^2,
    \label{eq::mzdmd_final_opt_unknown_B}
\end{equation}
where $\bX_+ = [\bx_m, \ldots, \bx_2]$ and $\bX_- = [\bx_{m-1}, \ldots, \bx_1]$ are the stacked snapshots and $\bV_c = [\bv_{m-1}, \ldots, \bv_1]$ is the stacked input vectors.
In this case, we again observe the one-to-one correspondence between our approach and the DMDc method up to memory term $\Delta \tau^2 \widetilde{\bM}(\bA, \bn)$.

It is to be noted here that this analysis presumes that the control $\bg$ does not explicitly depend on the solution. 
Otherwise, as follows from~\eqref{eq::caushy_char}, the control should lead to a modified operator $\mathcal{L}$ in the Liouville equation~\eqref{eq::liouv_general}. 
To take into account the control explicitly, we can presume that $\calB \bv = \bg$. 
This is always valid at $t=0$ and approximately valid for small enough~$t$. 
Then, we count the control in the framework of the first-order optimal prediction, when the memory term is neglected~\cite{chorin2002optimal}. 
This approach is justified at least for small enough time.

In the next section, we demonstrate that the DMDc method can not accurately predict the averaged (expected) trajectory for the resolved variables if it is based on a single trajectory.
In contrast, the proposed optimal prediction with control method quite accurately tracks the averaged trajectory from a single measurement for the resolved variables.

\section{Numerical experiments}
\label{sec:experiments}

In this section, we present an experimental evaluation of the proposed approach and compare it with the Monte Caro projection and DMDc. 
We consider both cases of known and unknown operator $\calB$.

\subsection{If the operator $\calB$ is known}
\label{sec::exp_known_B}
To solve the optimization problem~\eqref{eq::mzdmd_final_opt}, we use the Adam optimizer~\cite{kingma2014adam} and the automatic differentiation technique~\cite{baydin2018automatic} built in the JAX framework~\cite{jax2018github}.
The particular number of iterations for the optimizer and the learning rates are reported for each experiment separately.
We select the number of iterations such that the objective function does not change essentially.
To generate the ground-truth dynamics via the Monte Carlo (MC) projection, we average the trajectories obtained with 100 initialization of the unresolved variables.
At the same time, to generate the averaged dynamics from the computed operator $\bA_{opc}$, we use 1000 random vectors to compute the correction term in~\eqref{eq::opc_generation}.  
Next, we consider two test problems and evaluate the proposed method for the constant and damped control functions. 

\paragraph{Test problem 1.}
The first test problem we use to generate the data snapshots is the following
\begin{equation}
    \begin{split}
    \dot{y}_1 &= y_2+g_1\\
\dot{y}_2 &= -y_1(1 + y_3^2)+g_2\\
\dot{y}_3 &= y_4+g_3\\
\dot{y}_4 &= -y_3(1 + y_1^2)+g_4.\\
\end{split}
    \label{eq::test_problem1}
\end{equation}

The original dynamical system without the control is the Hamiltonian and corresponds to the coupled oscillators.
We use the standard normal distribution $\mathcal{N}(0, \mathbf{I})$ for the initialization of the unresolved variables $(y_3, y_4)$ and the initialization of the resolved variables is $(y_1(0), y_2(0)) = (1, 0)$.  

In these experiments, we use the following control functions to illustrate the robustness of the presented approach not only to the original dynamical system but also to the control function:
\begin{itemize}
    \item Constant control function $\mathbf{g}=\mathbf{g}_c \equiv \bc$, where $\bc$ is a constant vector.
    \item Damped control function $\mathbf{g}=\mathbf{g}_d(\by) \equiv -k \by.
$\end{itemize}

\begin{figure}[!ht]
    \centering
    \begin{subfigure}[t]{0.3\textwidth}
    \includegraphics[width=\textwidth]{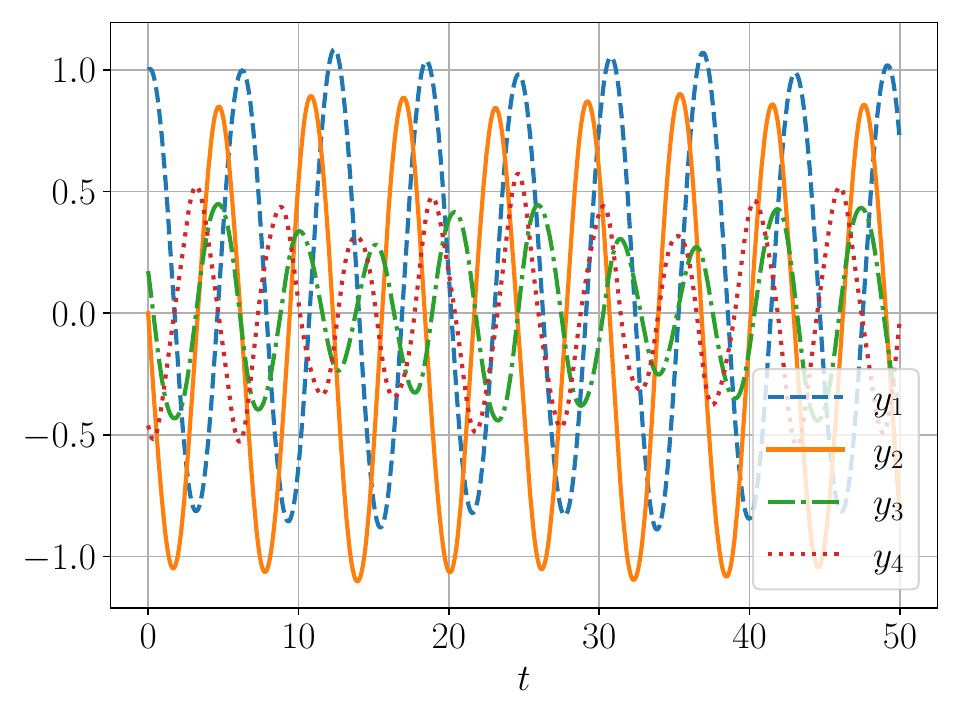}
    \caption{Measured trajectories}
    \label{fig::single_traj_const_test1}
    \end{subfigure}
    ~
    \begin{subfigure}[t]{0.3\textwidth}
    \includegraphics[width=\textwidth]{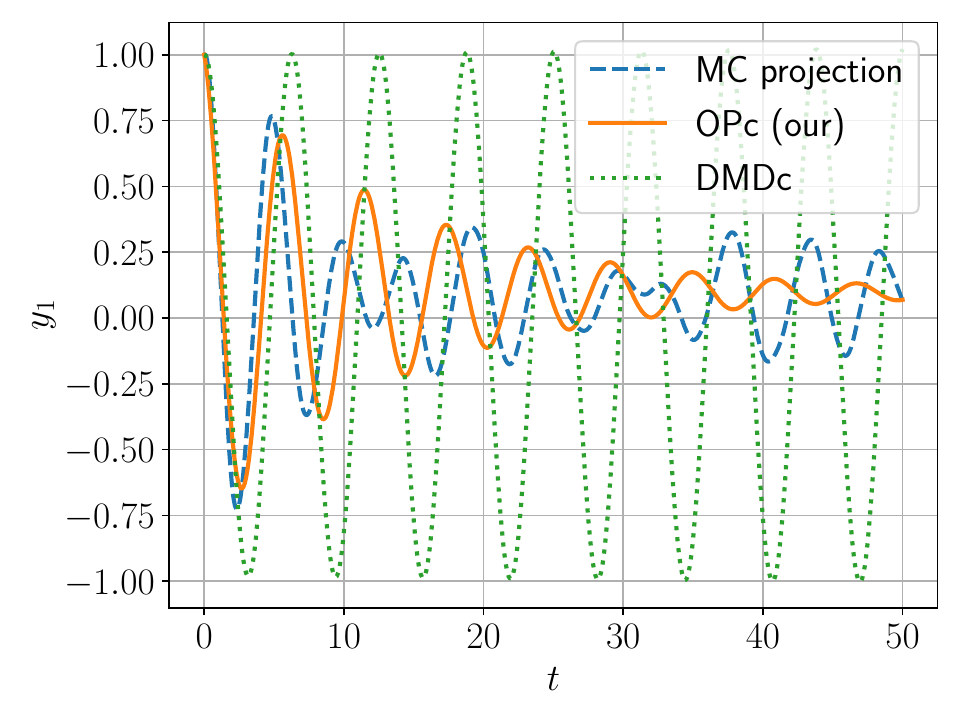}
    \caption{Comparison of averaged trajectories for $y_1$}
    \end{subfigure}
    ~
    \begin{subfigure}[t]{0.3\textwidth}
    \includegraphics[width=\textwidth]{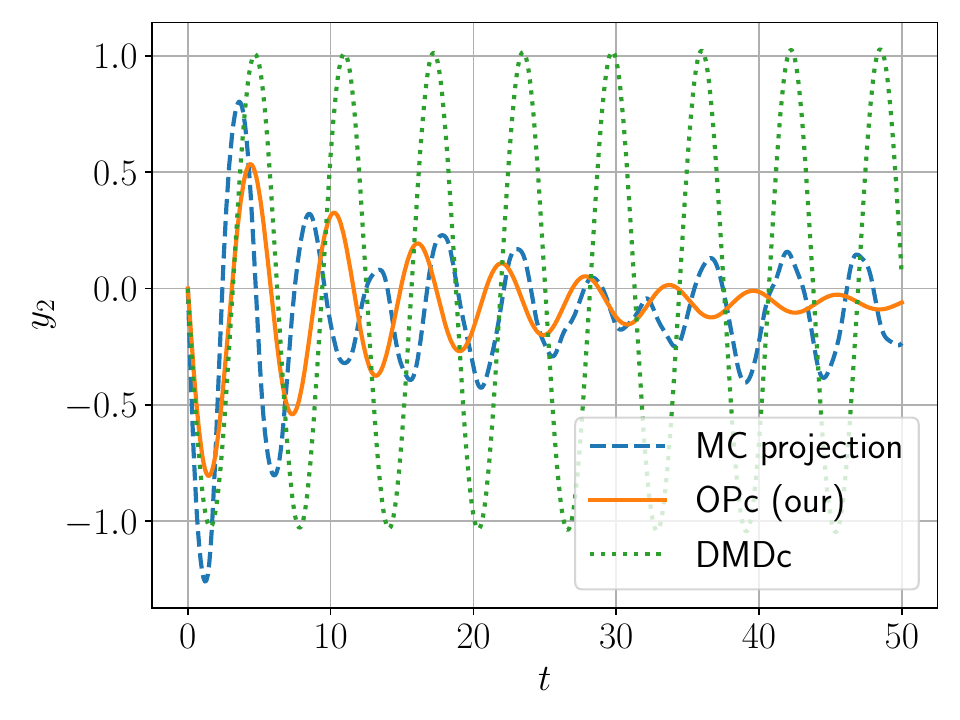}
    \caption{Comparison of averaged trajectories for $y_2$}
    \end{subfigure}
    \caption{Comparison of trajectories corresponding to the dynamical system~\eqref{eq::test_problem1} and the control function $\bg_c(\by) = [0.1, 0.1, -0.01, 0.01]$. (a): trajectories of resolved and unresolved variables corresponding to a single initialization. (b) and (c): dynamics obtained from the MC projection method and OPc method have a similar trend. This trend shows the stability of averaged trajectories corresponding to the resolved variables.
    A straightforward DMDc constructed from a single measurement shows regular trajectories, contradicting the ground-truth decay trajectories via the MC projection and OPc methods.
    }
    \label{fig::test_problem1_const}
\end{figure}

Let us start from constant control function $\bg_c(\by) = \bc$, where $\bc= [0.1, 0.1, -0.01, 0.01]$.
The vector $\bc$ values highlight the difference in the amplitudes of the resolved and unresolved variables.
Figure~\ref{fig::test_problem1_const} shows the regular dynamics derived with the DMDc method, which coincide with the measured trajectories, see Figure~\ref{fig::single_traj_const_test1}.
However, it does not approximate the averaged (expected) trajectory of the resolved variables given by the straightforward Monte Carlo (MC) projection.
This observation demonstrates a poor performance of DMDc based on the single measurement trajectory in the prediction of the expected solution. 
At the same time, the proposed OPc method provides the averaged trajectory of the resolved variables, which is similar to that generated via the MC projection.
In this experiment, we use 150 iterations of the Adam optimizer, and the learning rate equals $10^{-3}$.

    

\begin{figure}[!ht]
    \centering
    \begin{subfigure}[t]{0.3\textwidth}
    \includegraphics[width=\textwidth]{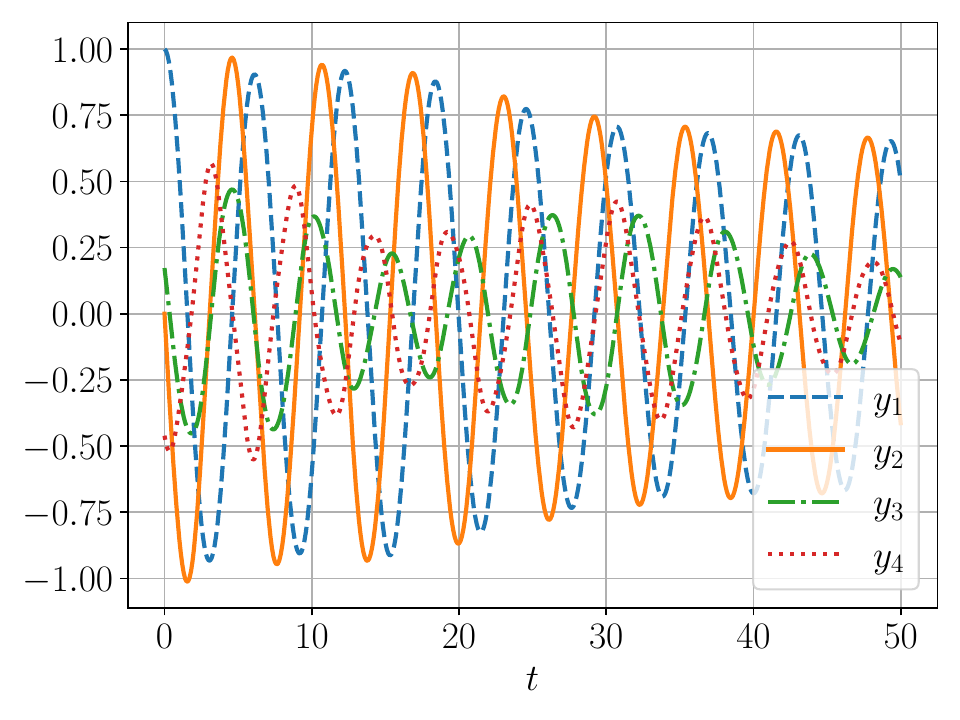}
    \caption{Measured trajectories}
    \end{subfigure}
    ~
    \begin{subfigure}[t]{0.3\textwidth}
    \includegraphics[width=\textwidth]{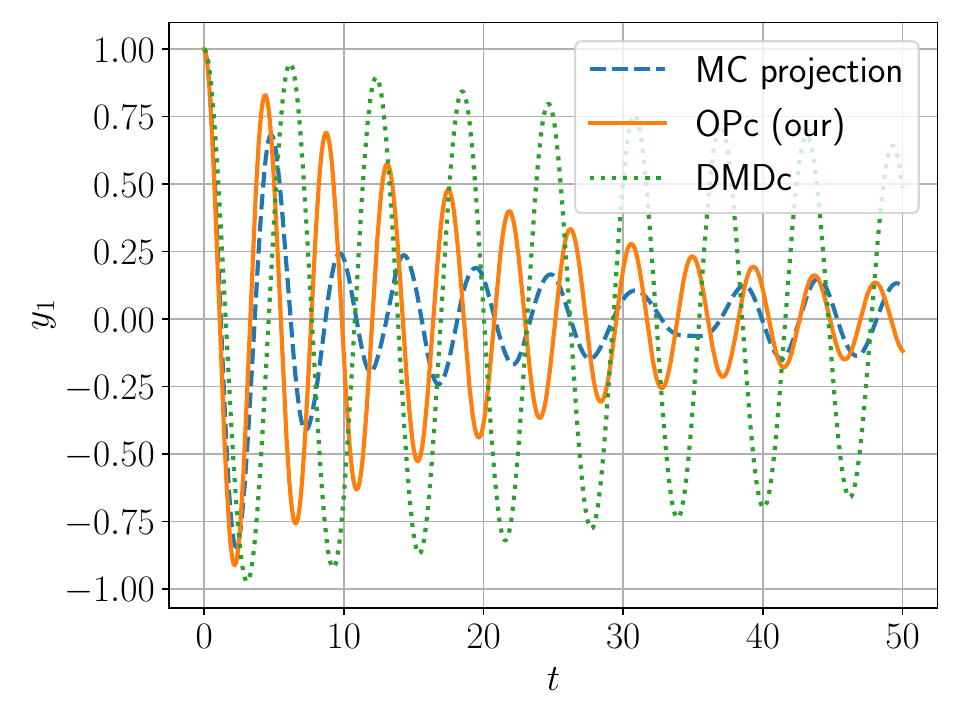}
    \caption{Comparison of averaged trajectories for $y_1$}
    \end{subfigure}
    ~
    \begin{subfigure}[t]{0.3\textwidth}
    \includegraphics[width=\textwidth]{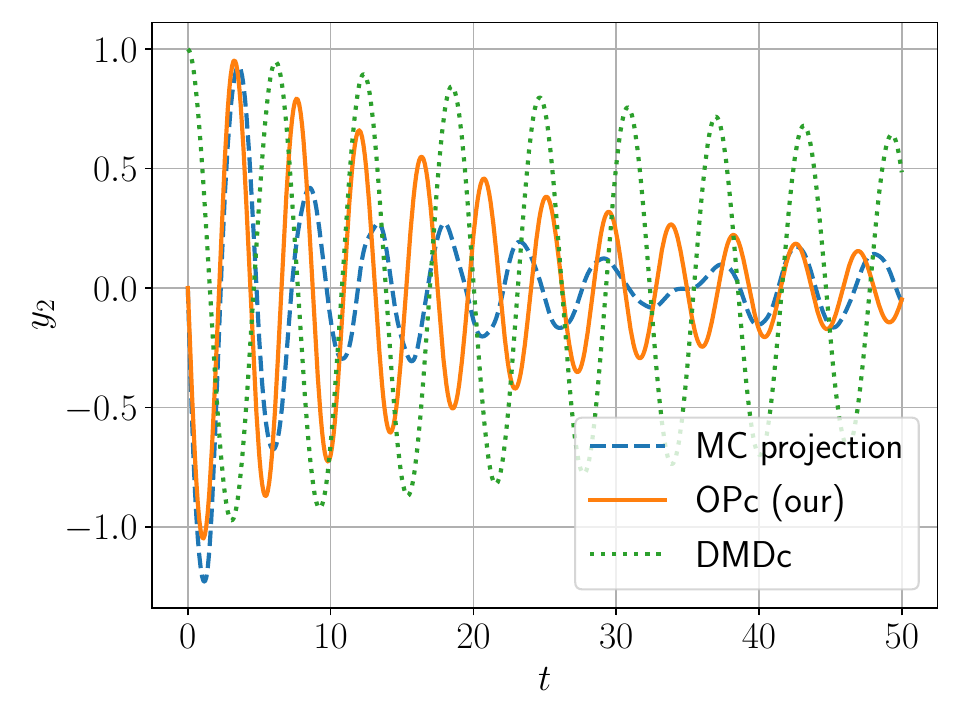}
    \caption{Comparison of averaged trajectories for $y_2$}
    \end{subfigure}
    \caption{Comparison of the trajectories corresponding to the dynamical system~\eqref{eq::test_problem1} and the control function $\bg_d(\by) = -0.01 \by$.
    (a): the trajectories of resolved and unresolved variables generated after a single initialization. (b) and (c): the proposed optimal prediction with control (OPc) method captures the same dynamics trend as the ground-truth solution computed via Monte Carlo (MC) projection. DMDc shows a significantly slower decay of the trajectories and does not provide an accurate approximation of the ground-truth averaged trajectories.}
    \label{fig::test_problem1_damped}
\end{figure}

In the next part the same dynamical system~\eqref{eq::test_problem1} is considered but with the damped control function $\bg_d(\by) = -0.01 \by$.
In case $\bg=\bg(\by)$, for certain $\calB \bv=\bg$ only at $t=0$. 
As such, if we set $\calB \bv=\bg$, then the operator $\bA_{opc}$ should slightly adjust to compensate for this inconsistency.
In this experiment, we use 200 iterations of the Adam optimizer and the learning rate is equal to $2\cdot 10^{-2}$.
Figure~\ref{fig::test_problem1_damped} demonstrates the reconstructed dynamics with the straightforward projection via the Monte Carlo simulations and with the proposed OPc method.
We also provide the auxiliary trajectories reconstructed by the DMDc method to illustrate its poor performance.
These plots confirm that the proposed approach can recover the correct trend even if only a single trajectory measurement is used.

\paragraph{Test problem 2.}
The second test problem is based on a modification of~\eqref{eq::test_problem1} and incorporates different scales for resolved and unresolved variables through the parameter $\varepsilon$:
\begin{equation}
\begin{split}
    \dot{y}_1 &= y_2+g_1\\
\dot{y}_2 &= -y_1 (1 + \varepsilon y_3^2)+g_2\\
\dot{y}_3 &= \varepsilon y_4+g_3\\
\dot{y}_4 &= -\varepsilon y_3 (1 + y_1^2)+g_4.\\
\end{split}
    \label{eq::test_problem_2}
\end{equation}
Here, we use the same control functions $\bg$ as those used in the previous test case.

In our tests, we use $\varepsilon = 10$ to simulate significantly different scales for the resolved and unresolved variables.
The corresponding initialization for the unresolved variables is generated from the normal distribution $\mathcal{N}(0, \sigma^2 \bI)$, where $\sigma = 1 / \sqrt{\varepsilon}$ according to~\cite{curtis2021dynamic}.
The initialization of the resolved variables is $(y_1(0), y_2(0)) = (1, 0)$ similar to the test problem 1. 

\begin{figure}[!ht]
    \centering
    \begin{subfigure}[t]{0.3\textwidth}
        \includegraphics[width=\textwidth]{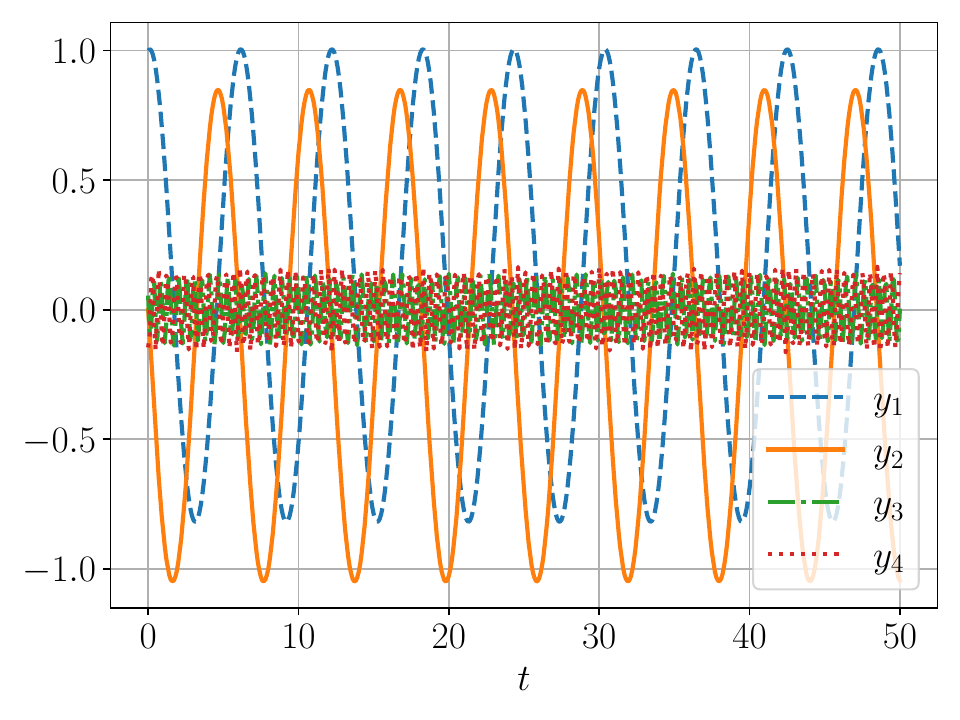}
        \caption{Measured trajectories}
        \label{fig::single_traj_const_test2}
    \end{subfigure}
    ~
    \begin{subfigure}[t]{0.3\textwidth}
    \includegraphics[width=\textwidth]{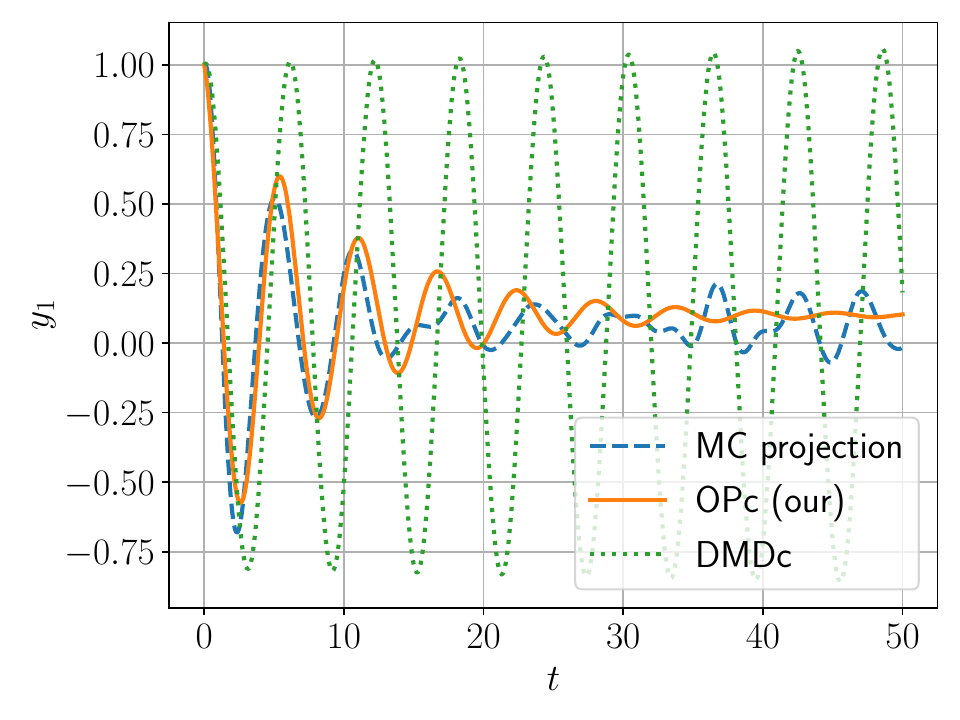}
    \caption{Comparison of averaged trajectories for $y_1$}
    \label{fig::cmp_test2_y1}
    \end{subfigure}
    ~
    \begin{subfigure}[t]{0.3\textwidth}
    \includegraphics[width=\textwidth]{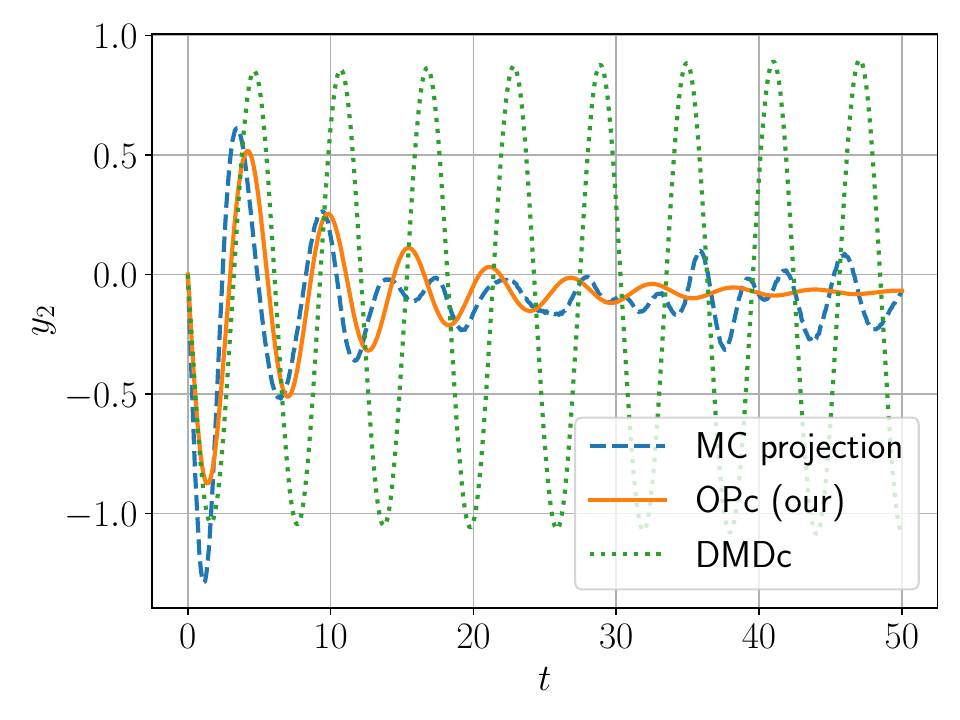}
    \caption{Comparison of averaged trajectories for $y_2$}
    \label{fig::cmp_test2_y2}
    \end{subfigure}
    \caption{Comparison of trajectories for the constant control function $\bg(\by) = \mathbf{c}$ and dynamical system~\eqref{eq::test_problem_2}. (a): trajectories of  resolved and unresolved variables with a single initialization. 
    (b) and (c): comparison of trajectories computed with MC projection, DMDc, and OPc methods for the resolved variables. 
    The OPc method captures the trend of ground-truth averaged trajectories given by MC projection.
    The DMDc method gives stable trajectories that differ from the observed decay behavior of the averaged trajectories from both OPc and MCc projection methods.
    }
    \label{fig::test2_constant_control}
\end{figure}

Let us start with the constant external control used for test problem 1.
In this setup, we run 150 iterations of the Adam optimizer with the learning rate equal to $10^{-3}$.
Figure~\ref{fig::test2_constant_control} shows that the OPc method recovers the averaged trajectory consistently with the ground-truth trajectory obtained from the Monte Carlo (MC) projection.
The DMDc method again shows a poor performance and disables to reconstruct the expected trajectory from a single measurement.
Note that the frequencies of the unresolved variables in the single measurement are much higher than those of the resolved variables; cf Figures~\ref{fig::single_traj_const_test2} and~\ref{fig::single_traj_const_test1}.

\begin{figure}[!ht]
    \centering
    \begin{subfigure}[t]{0.3\textwidth}
        \includegraphics[width=\textwidth]{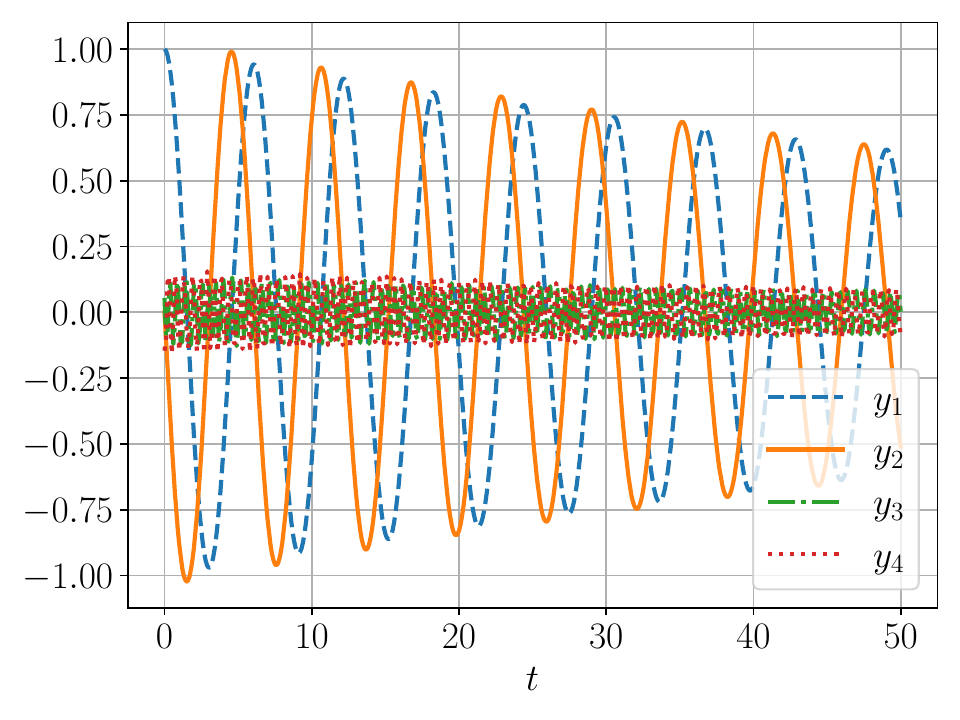}
        \caption{Measured trajectories}
        \label{fig::single_traj_damped_test2}
    \end{subfigure}
    ~
    \begin{subfigure}[t]{0.3\textwidth}
    \includegraphics[width=\textwidth]{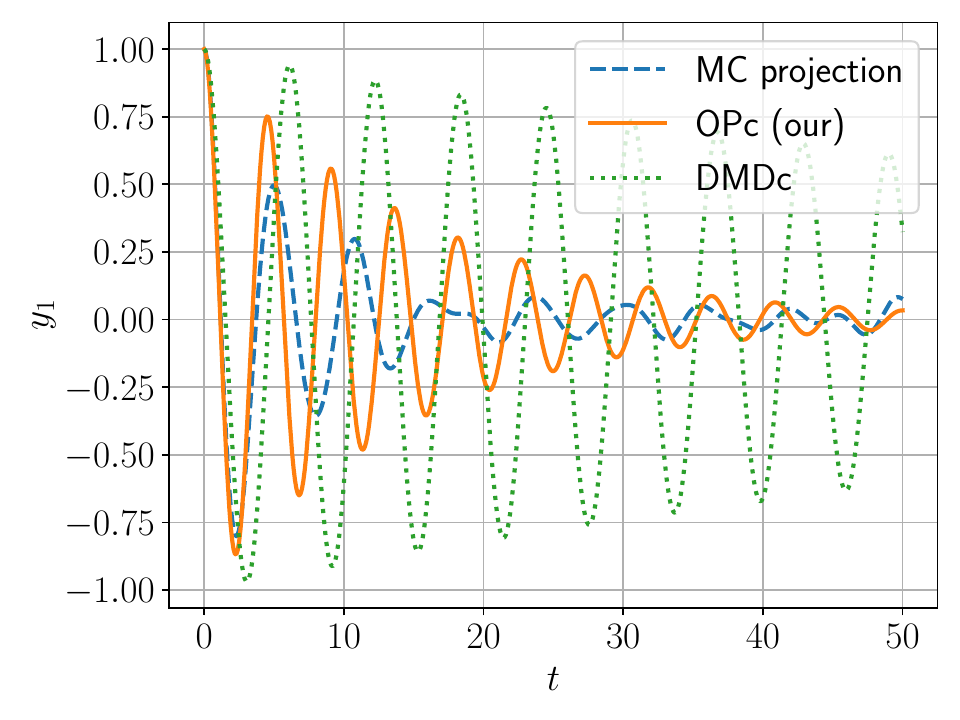}
    \caption{Comparison of averaged trajectories for $y_1$}
    \label{fig::cmp_test2_damped_y1}
    \end{subfigure}
    ~
    \begin{subfigure}[t]{0.3\textwidth}
    \includegraphics[width=\textwidth]{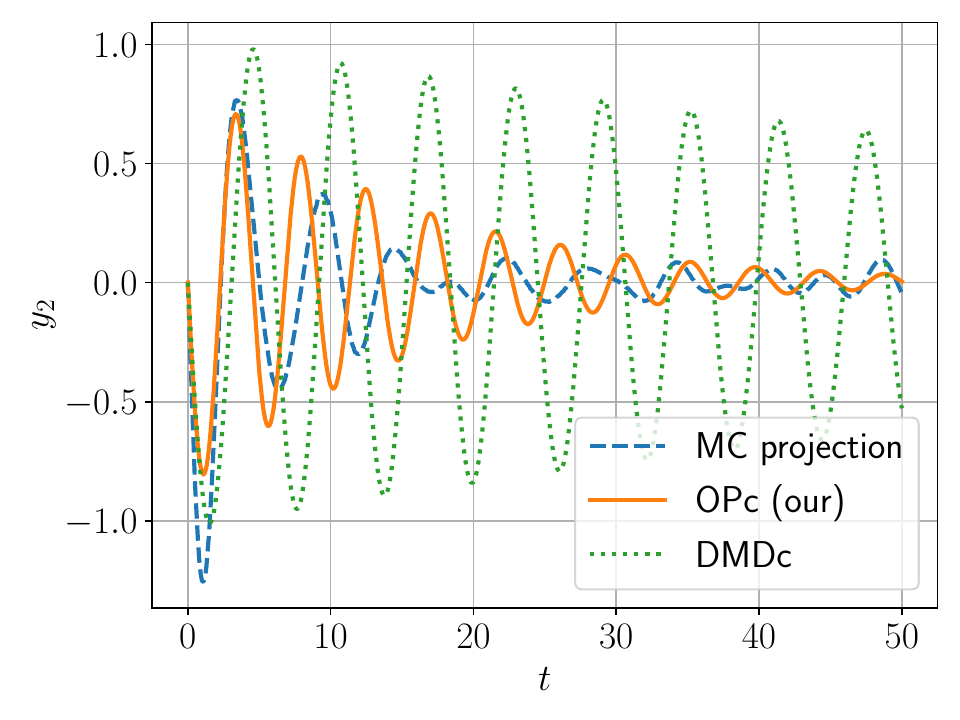}
    \caption{Comparison of averaged trajectories for $y_2$}
    \label{fig::cmp_test2_damped_y2}
    \end{subfigure}
    \caption{Comparison of trajectories for the damped control function $\bg_d(\by) = -0.01 \by$ and dynamical system~\eqref{eq::test_problem_2}. (a): trajectories of  resolved and unresolved variables with a single initialization. 
    (b) and (c): comparison of trajectories computed with MC projection, OPc, and DMDc methods for the resolved variables. 
    The OPc method captures the trend of the ground-truth averaged trajectories given by MC projection.
    DMDc gives slow decay trajectories that significantly differ from the expected trajectories.
    }
    \label{fig::test2_damped_control}
\end{figure}

The second tested control function is the damped control $\bg_d(\by) = -0.01\by$ used in the previous experiments.
Similarly to the previous setup, we run 300 iterations of the Adam optimizer, and the learning rate equals $2\cdot 10^{-2}$.
Figure~\ref{fig::single_traj_damped_test2} shows the effect of the damped control function on the trajectories obtained from the single initialization of the unresolved variables.
At the same time, Figures~\ref{fig::cmp_test2_damped_y1} and \ref{fig::cmp_test2_damped_y2} demonstrate that the OPc method accurately tracks the trend of the averaged trajectories obtained from the MC projection method.
One can also observe in Figures~\ref{fig::cmp_test2_damped_y1} and~\ref{fig::cmp_test2_damped_y2} that DMDc gives very slowly convergent trajectories which do not align with the behavior of the averaged trajectories.  

\subsection{If the operator $\calB$ is unknown}
To illustrate this setup, we take the test problem 1 (see~\eqref{eq::test_problem1}) and consider two possible cases:
\begin{enumerate}
    \item[1)] available input vectors $\bv_n$ have elements corresponding to both resolved and unresolved variables.
    We use complete state vector $\by_n$ as $\bv_n$ for simplicity.
    In this case, we generate the data using the matrix $\bB_1 = \diag{[-10^{-2}, -10^{-2}, -10^{-2}, -5 \cdot 10^{-2}]}$ and the control function in~\eqref{eq::caushy_char} looks like $\mathbf{g}(\by) = \bB_1\by$.
    \item[2)] available input vectors $\bv_n$ have elements corresponding to \emph{only} the resolved variables, and the control function $\mathbf{g}$ transforms these vectors to the control vectors affecting both the resolved and unresolved variables.
    As we use only the resolved variables $\bx_n$ as the input vector $\bv_n$, to generate the measured trajectories, we use the matrix $\bB_2 = \begin{bmatrix}
        -0.1 & 1 \\ -1 & 0 \\ -1 & 0.1 \\ -0.1 & -0.1
    \end{bmatrix}$ and, therefore, the control function in~\eqref{eq::caushy_char} looks like $\mathbf{g}(\by) = [\bB_2, \mathbf{0}]\by$.
\end{enumerate}
In both cases, we approximately set $\bB_{opc}\bv_n=\bg$.
However, since the operator $\bB$ is also tuned, the resulting $\bB_{opc}$ aims to approximate the ground-truth operator $\calB$ and correct the resulting control vectors $\bB_{opc}\bv_n$.
Thus, we expect the averaged trajectories to be reconstructed more accurately than in the previous case, considered in Section~\ref{sec::exp_known_B}.
After the measured trajectories are generated, OPc searches for the locally optimal matrices $\bA_{opc}$ and $\bB_{opc}$ based on the measured trajectories according to the optimization problem~\eqref{eq::mzdmd_final_opt_unknown_B}.
The target matrix $\bB_{opc}$ is in a general form, and its dimension coincides with the dimension of the space of resolved variables.

\begin{figure}[!ht]
    \centering
    \begin{subfigure}[t]{0.3\textwidth}
        \includegraphics[width=\linewidth]{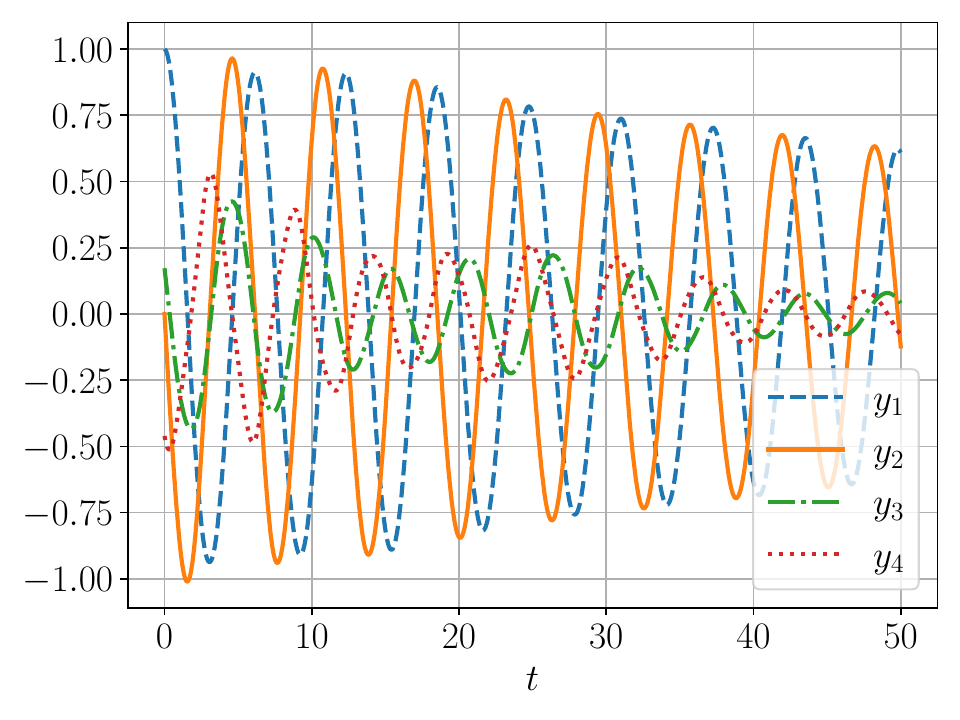}
        \caption{Measured trajectories}
        \label{fig::single_measurement_diag_B}
    \end{subfigure}
    ~
    \begin{subfigure}[t]{0.3\textwidth}
        \includegraphics[width=\linewidth]{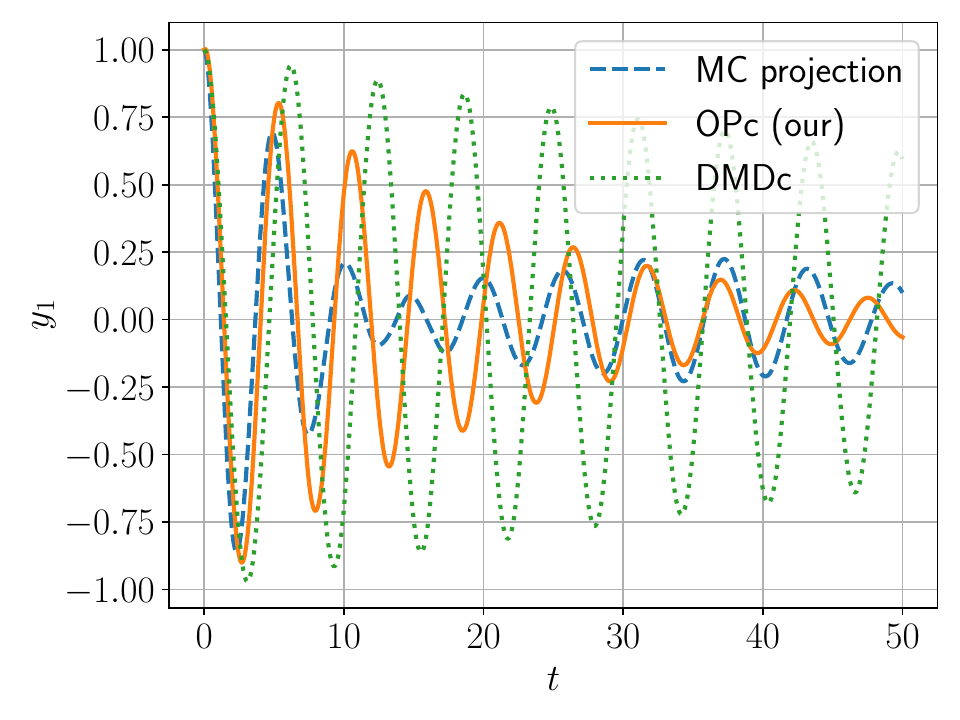}
        \caption{Comparison of averaged trajectories for $y_1$}
    \end{subfigure}
    ~
    \begin{subfigure}[t]{0.3\textwidth}
        \includegraphics[width=\linewidth]{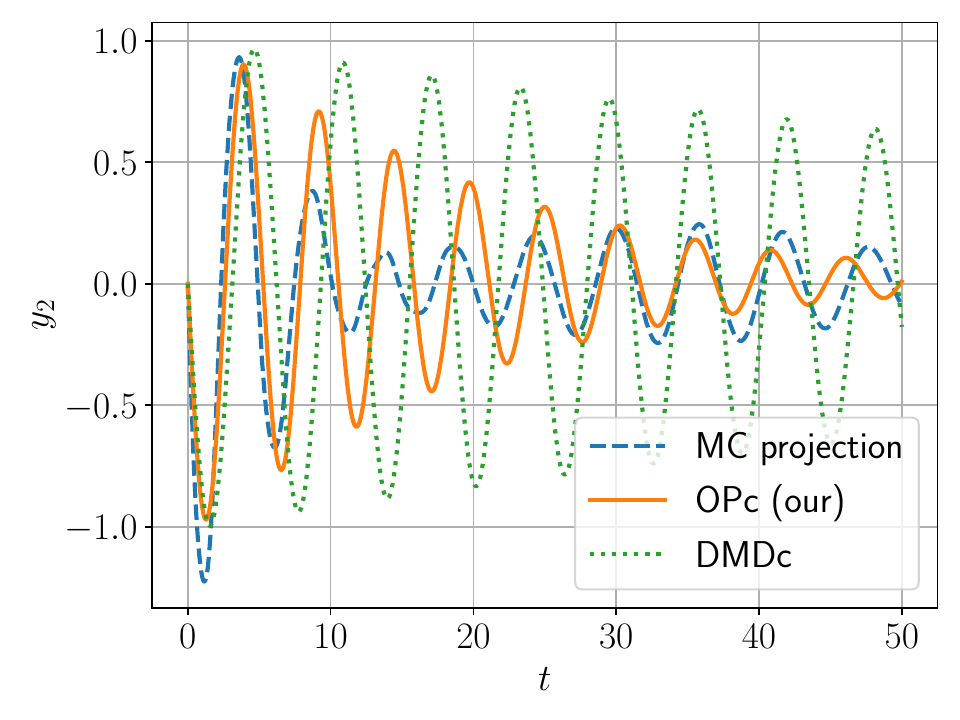}
        \caption{Comparison of averaged trajectories for $y_2$}
    \end{subfigure}
    \caption{Comparison of measured and averaged trajectories. (a): measurement trajectories show the decay of resolved and unresolved variables but with slightly different rates due to the values in the diagonal of matrix~$\bB_1$. (b) and (c): comparison of reconstructed average trajectories with the expected trajectories confirms that the proposed OPc method can recover matrix $\bB_{opc}$. The latter corresponds to the given input vectors affecting resolved and unresolved variables.}
    \label{fig::test1_diag_B}
\end{figure}

Let us start with the first case, where the input vectors coincide with the complete state vectors.
We use 300 iterations of Adam optimizer and the learning rate equal to $10^{-3}$.
Figure~\ref{fig::test1_diag_B} shows that the proposed OPc method recovers both $\bA_{opc}$ and $\bB_{opc}$ such that the recovered trajectories are close to the ground-truth trajectories computed via the MC projection.
We also observe that the DMDc trajectories are quite close to the measured trajectories but they significantly differ from the averaged trajectories.
Thus, we demonstrate that the OPc method can reconstruct not only the matrix $\bA_{opc}$, which approximates the projected operator $\mathcal{L}$ in~\eqref{eq::liouv_general}, but also matrix~$\bB_{opc}$.

\begin{figure}[!ht]
    \centering
    \begin{subfigure}[t]{0.3\textwidth}
        \includegraphics[width=\linewidth]{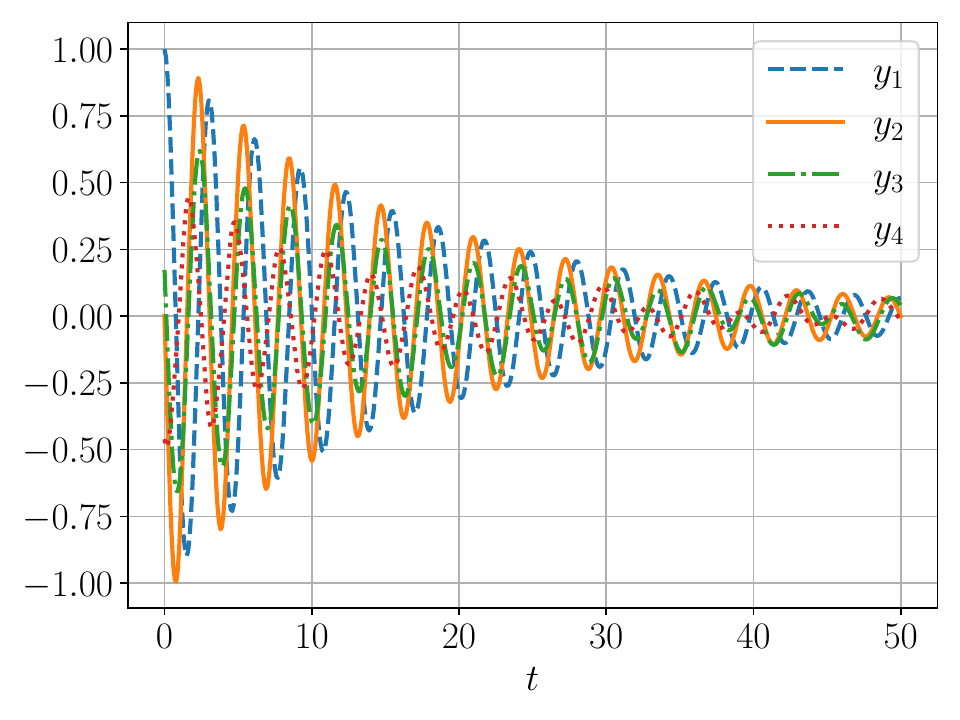}
        \caption{Measured trajectories}
    \end{subfigure}
    ~
    \begin{subfigure}[t]{0.3\textwidth}
        \includegraphics[width=\linewidth]{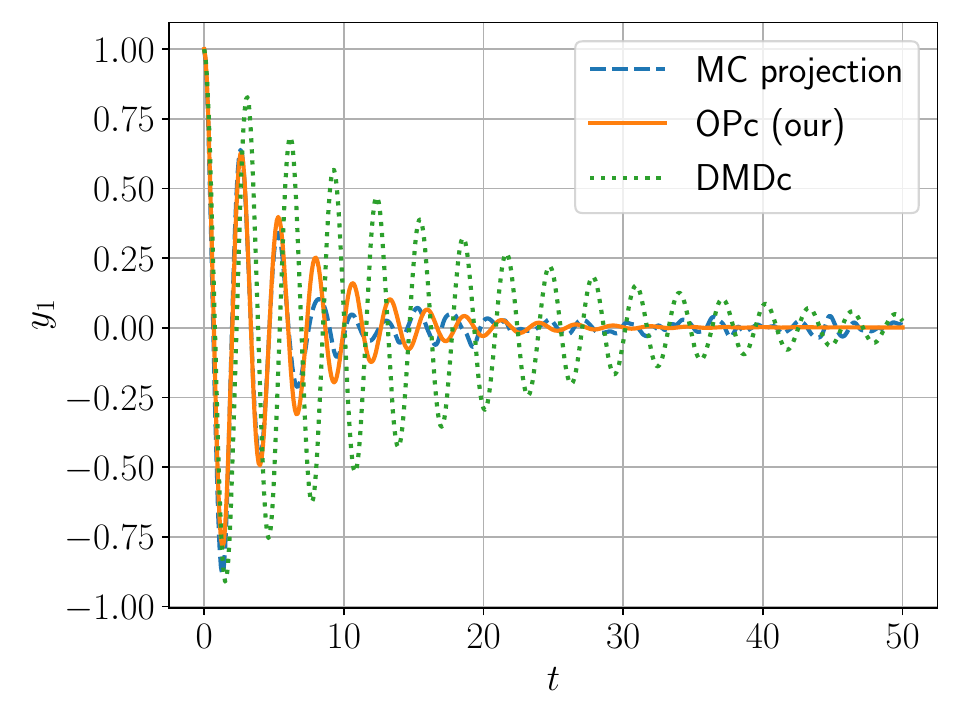}
        \caption{Comparison of the averaged trajectories for $y_1$}
    \end{subfigure}
    ~
    \begin{subfigure}[t]{0.3\textwidth}
        \includegraphics[width=\linewidth]{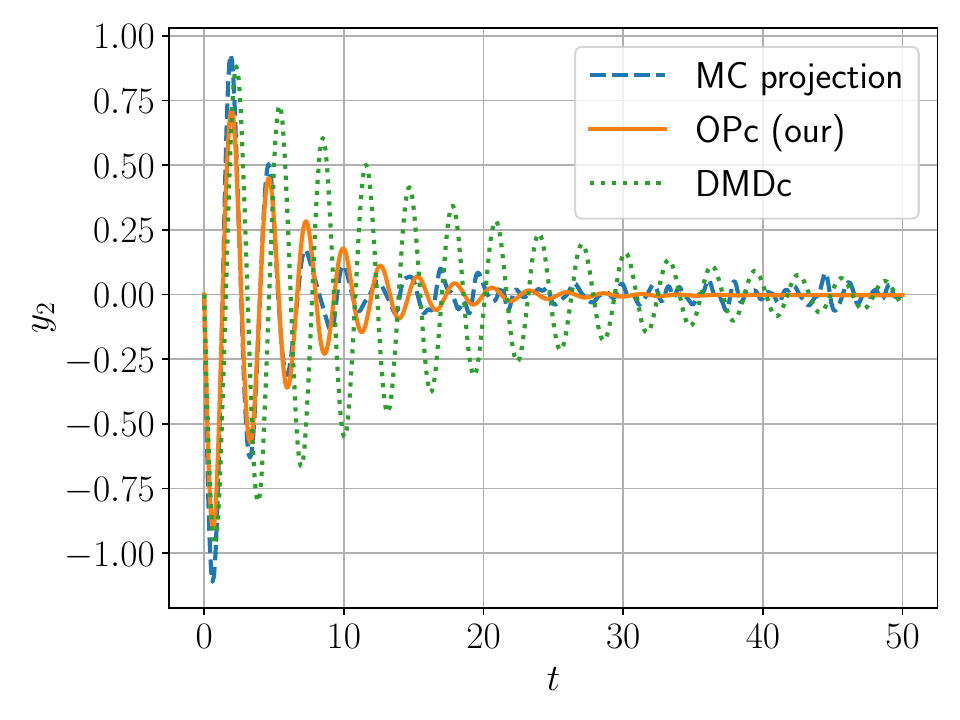}
        \caption{Comparison of the averaged trajectories for $y_2$}
    \end{subfigure}
    \caption{
    Comparison of measured and averaged trajectories. 
    (a): measurement trajectories show the decay of both resolved and unresolved variables, which become very close by the end of the trajectories. (b) and (c): the comparison of the reconstructed average trajectories with the expected trajectories demonstrates that the proposed OPc method gives the accurate reconstruction of the unknown matrices $\bA_{opc}$ and $\bB_{opc}$. The latter corresponds to the given input vectors affecting only resolved variables.
    }
    \label{fig::test1_rect_B}
\end{figure}

The second case corresponds to the reduced input vectors to the resolved variables only.
In this experiment, we use 400 iterations of the Adam optimizer and the learning rate equal to $10^{-3}$.
Figure~\ref{fig::test1_rect_B} shows that the measured trajectories decay for both the resolved and unresolved variables, and their amplitudes become very close for $t > 35$.
Thus, we see the difference with the previous case, cf.~Figure~\ref{fig::single_measurement_diag_B}.
We can also see from Figure~\ref{fig::test1_rect_B} that the OPc method provides a very accurate reconstruction of the averaged trajectories that are almost indistinguishable from the MC projection results for $t < 7$.
At the same time, the DMDc method again provides slow decayed trajectories that are not aligned with the MC projection trajectories.

\section{Discussion}

In the previous section, we have focused on the accuracy of the averaged trajectories reconstruction and the details of the experimental setup.
However, the accuracy of reconstruction is not enough to be practically important since one can obtain the averaged trajectories using Monte Carlo simulations.
Therefore, we present the runtime comparison of the considered approaches for the benchmarks from the previous section in Table~\ref{tab::runtime_cmp}.
We present the runtime only for the case of the known operator $\calB$ and note that the runtime for the search matrices $\bA_{opc}$ and $\bB_{opc}$ is very close to the reported.
This table shows that our approach generates the averaged trajectories from the single measurement much faster than the Monte Carlo projection.
Note that we report only the pure runtime and ignore timing spent on compiling the Python code with the JIT technique built in the JAX framework.

\begin{table}[!ht]
    \centering
    \caption{Summary of comparing the method OPc with straightforward Monte Carlo projection for the considered benchmarks. The runtime is reported in the table in seconds.}
    \begin{tabular}{ccccc}
    \toprule
        \multirow{2}{*}{Methods} &  \multicolumn{2}{c}{Test problem 1} & \multicolumn{2}{c}{Test problem 2}\\
        & $g_d$ & $g_c$ & $g_d$ & $g_c$ \\
        \midrule
        MC projection & $11.64$ & $13.18$ & $54.89$ & $20.47$ \\
        OPc & $0.97$ & $0.92$ & $1.03$ & $0.97$ \\
         \bottomrule
    \end{tabular}
    \label{tab::runtime_cmp}
\end{table}

In addition to the runtime comparison, we discuss the place of the OPc method among the other data-driven methods for the analysis of the dynamical systems.
We summarize the main use cases in~Table~\ref{tab::context}.
This table illustrates that OPc complements the optimal prediction (OP) and DMD-type techniques.
Only the OPc method is applicable for analyzing the dynamical systems with control and unresolved variables in the state vector.
Thus, the proposed OPc method is theoretically and practically important for predicting the average trajectories of the dynamical systems with control under uncertainties.

\begin{table}[!ht]
    \centering
    \caption{The summary of use cases for the proposed OPc method and the related ones. The key questions to decide what method is appropriate are: do we have only the resolved variables, and do we have a control function in our model of the dynamical system.}
    \begin{tabular}{cc|cc}
      &  & \multicolumn{2}{c}{Control function} \\
     &  & Yes & No \\
    \hline
  \multirow{3}{*}{\rotatebox[origin=c]{90}{\parbox{1.5cm}{\centering Resolved only}}} & Yes & DMDc & DMD \\
  & & & \\
   & No & OPc & OP \\
    \end{tabular}
    \label{tab::context}
\end{table}

\section{Conclusion and future work}
\label{sec:conclusions}

This study presents the OPc method, which extends the optimal prediction approach to the dynamical systems with control functions.
The main feature of the considered dynamical systems is that the state vector depends on unresolved variables typically modeled via random initializations.
Our data-driven method reconstructs the averaged trajectories of the resolved variables from a single measurement based on a random initialization of unresolved variables.
We illustrate the accuracy of the prediction of averaged trajectories for different test dynamical systems equipped with damped and constant control functions.
The case of the unknown control operator is also addressed by the proposed OPc method, and the averaged trajectories are also reconstructed accurately enough.   
In addition to the reconstruction accuracy, we demonstrate that the proposed OPc method is much faster than the straightforward Monte Carlo projection.
On the top, the OPc method provides the locally optimal linear operator whose spectrum can be used to evaluate the stability of the averaged trajectories for the resolved variables. 

Further theoretical analysis of the optimization problem and derivation of the convergence rate for an optimizer are interesting research directions that can be studied in the future work.
The proposed OPc method can also be extended to more complicated non-linear regimes with control.
For example, the dynamics of the control vector may be governed by the separated set of PDEs.
To tackle such a problem, one may use the PDE-constrained optimization framework.


\bibliographystyle{unsrt}
\bibliography{lib}

\end{document}

%% file: math_symbols.tex
\newcommand{\bx}{\mathbf{x}}
\newcommand{\bu}{\mathbf{u}}
\newcommand{\bw}{\mathbf{w}}

\newcommand{\bs}{\mathbf{s}}

\newcommand{\bz}{\mathbf{z}}
\newcommand{\by}{\mathbf{y}}
\newcommand{\bc}{\mathbf{c}}

\newcommand{\bg}{\mathbf{g}}
\newcommand{\bv}{\mathbf{v}}
\newcommand{\bn}{\mathbf{n}}
\newcommand{\bA}{\mathbf{A}}
\newcommand{\bB}{\mathbf{B}}
\newcommand{\bI}{\mathbf{I}}

\newcommand{\bM}{\mathbf{M}}
\newcommand{\bN}{\mathbf{N}}
\newcommand{\bF}{\mathbf{F}}

\newcommand{\bY}{\mathbf{Y}}
\newcommand{\bV}{\mathbf{V}}

\newcommand{\bLambda}{\boldsymbol{\Lambda}}

\newcommand{\bX}{\mathbf{X}}
\newcommand{\bbR}{\mathbb{R}}

\newcommand{\calB}{\mathcal{B}}

\newcommand{\calP}{\mathcal{P}}
\newcommand{\calM}{\mathcal{M}}
\newcommand{\calQ}{\mathcal{Q}}

\newcommand{\calI}{\mathcal{I}}
\newcommand{\calL}{\mathcal{L}}

\newcommand{\diag}[1]{\mathrm{diag}(#1)}

\DeclareMathOperator*{\argmin}{\arg\,\min}